\documentclass[preprint,12pt]{elsarticle}

\usepackage{amssymb}
\usepackage{amsmath}

\begin{document}

\begin{frontmatter}



\title{Discrete-Time Fractional-Order PID Controller: Definition, Tuning, Digital Realization and Experimental Results\footnote{Please be advised that this a preprint submitted to a journal for possible publication and it may be subjected to some changes during the review process.}}


\author{Farshad Merrikh-Bayat$^{a,\dag}$, Seyedeh-Nafiseh Mirebrahimi$^a$, Mohammad-Reza Khalili$^a$}

\address{$^a$Department of Electrical and Computer Engineering, University of Zanjan, Zanjan,
    Iran, e-mail: {\{}f.bayat, n.mirebrahimi{\}}@znu.ac.ir,
    adasbios2@gmail.com\\$^\dag$ Corresponding Author, Tel: +98 (241) 5154061}

\begin{abstract}
In some of the complicated control problems we have to use the
controllers that apply nonlocal operators to the error signal to
generate the control. Currently, the most famous controller with
nonlocal operators is the fractional-order PID (FOPID). Commonly,
after tuning the parameters of FOPID controller, its transfer
function is discretized (for realization purposes) using the
so-called generating function. This discretization is the origin
of some errors and unexpected results in feedback systems. It may
even happen that the controller obtained by discretizing a FOPID
controller works worse than a directly-tuned discrete-time
classical PID controller. Moreover, FOPID controllers cannot
directly be applied to the processes modeled by, e.g., the ARMA or
ARMAX model. The aim of this paper is to propose a discrete-time
version of the FOPID controller and discuss on its properties and
applications. Similar to the FOPID controller, the proposed
structure applies nonlocal operators (with adjustable memory
length) to the error signal. Two methods for tuning the parameters
of the proposed controller are developed and it is shown that the
proposed controller has the capacity of solving complicated
control problems with a high performance.
\end{abstract}


\end{frontmatter}


\section{Introduction}
During the past seven decades PID controllers have been
successfully used in a wide variety of industrial applications
\cite{astrom2,visioli}. Currently, various continuous and
discrete-time versions of this type of controller are available,
which can be applied to the processes modeled by linear
differential or linear difference equations, respectively (which,
of course, correspond to rational transfer functions in the
complex $s$ or $z$ variable, respectively)
\cite{astrom2}-\cite{ogata}. Successful applications of PID
controllers to nonlinear processes can also be found in the
literature \cite{chang,lin}.

According to the high achievement and the simplicity of design and
implementation of PID controllers, many researchers tried to
enhance the performance of these controllers by innovating new
structures and tuning methods \cite{astrom1}-\cite{bucklaew}.
During the years, the main motivation behind all these attempts
had been the increasing need to develop easy-to-use and
high-performance controllers with the capacity of solving more
complicated control problems (for example, it is well understood
that classical PID controllers do not lead to satisfactory results
when the process is non-minimum phase, high order, non-linear or
has very large dead times \cite{astrom2}). One of these attempts
led to a new generation of PID controllers, known as the
fractional-order PID (FOPID) or $\mathrm{PI^\lambda D^\mu}$
controller, which was first proposed by Podlubny in 1999
\cite{podlubny}. In FOPID controllers the error and control,
respectively denoted as $e(t)$ and $u(t)$, are related from the
following equation
\begin{equation}\label{fopid1}
u(t)=k_p e(t)+k_iD_t^{-\lambda} e(t)+k_dD_t^\mu e(t),
\end{equation}
where $k_p,k_i,k_d\in\mathbb{R}$ and $\lambda,\mu\in\mathbb{R}^+$
are the parameters of controller to be tuned, and $D_t^{-\lambda}$
and $D_t^\mu$ are the fractional integral and differential
operator respectively, often defined by the Riemann-Liouville
definition as the following \cite{podlubny2}:
\begin{equation}\label{int}
D_t^{-\lambda} f(t)=\frac{1}{\Gamma(\lambda)} \int_0^t
\frac{f(\tau)}{(t-\tau) ^{1-\lambda}}\mathrm{d}\tau,
\end{equation}
\begin{equation}\label{diff}
D^\mu_t f(t)=\frac{1}{\Gamma(m-\mu)} \left(\frac{\mathrm{d}}
{\mathrm{d}t}\right)^m \int_0^t \frac{f(\tau)}{(t-\tau)
^{1+\mu-m}}\mathrm{d}\tau,
\end{equation}
where $m$ is a positive integer such that $m-1<\mu\le m$ and
$\Gamma(.)$ is the well known gamma function. Other definitions
for fractional differential operators can also be found in the
literature \cite{podlubny2}. Note that according to (\ref{int})
and (\ref{diff}), fractional integral and derivative are nonlocal
operators which apply the past values of $f$ to determine
$D_t^{-\lambda} f(t)$ and $D^\mu_t f(t)$. Hence, unlike classical
PID controllers, the derivative term of the FOPID controller is
actually a nonlocal operator acts on the error signal. Taking the
Laplace transform from both sides of (\ref{fopid1}) leads to the
following transfer function for FOPID controller \cite{podlubny}:
\begin{equation}\label{fopid_transfer}
C(s)=\frac{U(s)}{E(s)}=k_p+k_is^{-\lambda}+k_d s^\mu.
\end{equation}

Here it is worth to mention that unlike the classical PID
controllers, currently there is no direct definition available for
discrete-time FOPID controllers. However, for realization
purposes, it is common practice to first design the
(continuous-time) FOPID controller and then approximate it with a
discrete-time system. This approximation is often performed by
using the so-called \emph{generating function}
\cite{machado}-\cite{vinagre2}. In this technique the Laplace
variable $s$ in (\ref{fopid_transfer}) is substituted with a
certain function of $z$ (using, e.g., the Tustin transform) and
then the power series expansion (PSE) of the resulted expression
is obtained in terms of $z$. Finally, since any practical
discrete-time system must necessarily use a limited memory, the
resulted PSE is truncated. This approach suffers from many
drawbacks (see the discussions below), the source of all is the
unavoidable mismatch between the frequency response of continuous
and discrete transfer functions. The aim of this paper is to
propose a long-memory discrete-time PID (LDPID) controller which
removes the limitations of the above-mentioned discretization
scheme and still has the high performance of FOPID controllers.

There are many good reasons for defining and using LDPID
controllers. First of all, it is very common practice to model a
real-world continuous-time process by a (discrete) transfer
function in the $z$ variable. For example, such a model is
obtained when an unknown process is identified using ARMA or ARMAX
model. Obviously, in dealing with such a discrete-time model of
process it is more reasonable to directly design a discrete-time
controller as well (instead of designing a continuous-time
controller and then approximating it with a discrete-time one).
Another good reason for developing LDPID controllers is that even
when both the controller and process are continuous-time, the
controller is more likely to be realized using digital
microprocessors. As a classical fact, some unwanted effects may
occur in the feedback system when the continuous-time controller
is replaced with an approximate discrete-time controller (for
example, the phase margin is decreased according to the effect of
sample-and-hold and, in addition, the feedback system may become
unstable according to the effect of truncation).

Another important point in relation to the classical approximation
methods is that FOPID controllers have some features which are not
preserved after approximating them with discrete transfer
functions. For example, one important feature of every FOPID
controller is the so-called \emph{long memory principle}, which is
lost after approximating it with an integer-order transfer
function. Another property that is lost after approximation is the
optimality of controller. In fact, the optimal controller designed
for a certain continuous-time process will no longer be optimal
after approximating it with a discrete-time controller.

In addition to the above-mentioned points, it should also be noted
that the derivative term of the FOPID controller given in
(\ref{fopid_transfer}) cannot exactly be realized in practice
since it is a non-causal operator. In fact this term should always
be considered in series with a low-pass filter to be realizable
(more precisely, the practical FOPID controller is actually a six
degrees of freedom controller). The proposed LDPID controller will
remove this difficulty.

The rest of this paper is organized as the following. In Section 2
we introduce the proposed LDPID controller and develop two methods
for its tuning. Four simulated examples and an experimental study
are presented in Section 3. Finally, Section 4 concludes the
paper.

\section{Long-Memory Discrete-Time PID Controller}
\subsection{Formulation of the proposed LDPID controller}\label{sec_ldpid}
Consider the FOPID controller given in (\ref{fopid_transfer}). In
the following, without a considerable loss of generality, first we
develop a method for discretizing the derivative term of this
controller based on the prewarped Tustin method and then we extend
the results to the integrative term. Next, based on these results
we will propose a LDPID controller that works very similar to
FOPID controllers.

\begin{figure}[tb]
\begin{center}
\includegraphics[width=8cm]{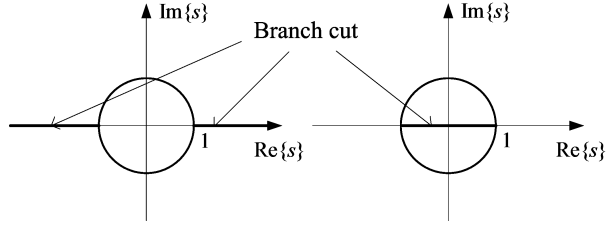}
\caption{Two possible branch cuts for the multi-valued function
given in (\ref{tus1}).}\label{fig_bc}
\end{center}
\end{figure}

Assuming that the sampling period of system is equal to $T$,
applying prewarped Tustin method leads to the following
approximation for the derivative term of FOPID controller:
\begin{equation}\label{tus1}
s^\mu=\left(\frac{\omega_c}{\tan(\omega_cT/2)}
\times\frac{1-z^{-1}}{1+z^{-1}} \right)^\mu = \alpha^\mu \times
\left(\frac{1-z^{-1}}{1+z^{-1}} \right)^\mu,
\end{equation}
where $\alpha\triangleq\omega_c/\tan(\omega_cT/2)$ and $\omega_c$
is the gain crossover frequency of the open-loop transfer
function. (Here we have applied Tustin method since it is more
accurate compared to other transforms such as backward and forward
difference. Section \ref{sec_general} makes a connection between
different possible transforms.) To proceed, we need to calculate
the PSE of the expression in the right hand-side of (\ref{tus1}),
which cannot be performed without determining its region of
convergence (ROC) in the complex $z$-plane. Note that this
expression is actually a multi-valued function of $z$ which has
two branch points at $z=1$ and $z=-1$ located on the unit circle.
Mathematically, the corresponding branch cut (BC) can be
considered either inside or outside the unit circle as shown in
Fig. \ref{fig_bc}. Clearly, the choice of BC affects the ROC and
consequently, causality of the resulted system. Similar to the
classical results \cite{oppenheim}, here considering the BC inside
the unit circle and the ROC as $|z|>1$ yields a causal system as
it is desired. Hence, assuming $w=z^{-1}$ and $|z|>1$ the PSE of
the expression in the right hand-side of (\ref{tus1}) is obtained
as the following:
\begin{equation}\label{expand1}
\alpha^\mu \times \left(\frac{1-w}{1+w} \right)^\mu= \alpha^\mu
\sum_{k=0}^\infty f_k(\mu) w^k,\quad |w|<1,
\end{equation}
where
\begin{equation}\label{weight}
f_k(\mu)=\frac{1}{k!}\times \frac{\mathrm{d}^k}
{\mathrm{d}w^k}\left(\frac{1-w}{1+w} \right)^\mu \bigg|_{w=0}.
\end{equation}
Substitution of $w=z^{-1}$ in (\ref{expand1}), and then the
resulted equation in (\ref{tus1}) yields the following PSE for the
fractional-order differentiator:
\begin{equation}\label{deriv_app}
s^\mu= \alpha^\mu \sum_{k=0}^\infty f_k(\mu) z^{-k},\quad |z|>1,
\end{equation}
where again the coefficients $f_k(\mu)$ are calculated from
(\ref{weight}). It can be shown (using Maple) that the first few
coefficients in (\ref{deriv_app}) are as the following:
$$f_0(\mu)=1,\quad f_1(\mu)=-2\mu,\quad f_2(\mu)=2\mu^2,$$
$$f_3(\mu)= -\frac{4}{3}\mu^3-\frac{2}{3}\mu,\quad f_4(\mu)= \frac{2}{3}\mu^4+\frac{4}{3}\mu^2,$$
$$f_5(\mu)=-\frac{4}{15}\mu^5-\frac{4}{3}\mu^3-\frac{2}{5}\mu,$$
$$\quad f_6(\mu)=\frac{4}{45}\mu^6+\frac{8}{9}\mu^4+\frac{46}{45}\mu^2,
~\ldots$$
Note that (\ref{deriv_app}) holds for both the positive and
negative values of $\mu$. Hence, one may try to expand the
integral term of (\ref{fopid_transfer}) in a similar manner and
arrive at an equation like (\ref{deriv_app}) in $\lambda$, but the
problem with such an expansion is that the resulted series does
not have infinite DC gain (considering the fact that any infinite
series must be truncated in practice), which is essential for
tracking the step command without steady-state error. In order to
find a series approximation for $s^{-\lambda}$ in terms of
$z^{-1}$ which has infinite DC gain, first we write it as
$s^{-\lambda}=(1/s)\times s^{1-\lambda}$ and then apply the
prewarped Tustin method to it. Applying this technique yields
\begin{equation}\label{int_app}
s^{-\lambda}= \alpha^{-\lambda} \frac{1+z^{-1}}{1-z^{-1}}
\sum_{k=0}^\infty f_k(1-\lambda) z^{-k},\quad |z|>1,
\end{equation}
where $f_k(1-\lambda)$ are again calculated from (\ref{weight}).

\begin{figure}[tb]
\begin{center}
\includegraphics[width=10cm]{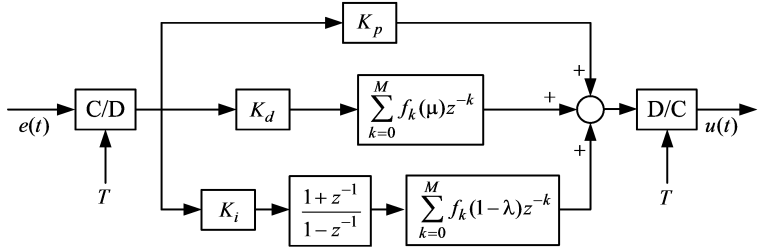}
\caption{Block diagram of the $M$th-order LDPID to be applied
instead of the continuous-time FOPID}\label{fig_dfopid_prac}
\end{center}
\end{figure}

Substitution of (\ref{deriv_app}) and (\ref{int_app}) in
(\ref{fopid_transfer}) results in the following formulation for
the LDPID controller
\begin{equation}\label{dfopid_ideal}
C_d(z)=K_p+K_d \sum_{k=0}^\infty f_k(\mu)z^{-k}+ K_i
\frac{1+z^{-1}}{1-z^{-1}} \sum_{k=0}^\infty f_k(1-\lambda)z^{-k},
\end{equation}
where
\begin{equation}\label{kp}
K_p=k_p,\quad K_d=k_d \alpha^\mu,\quad K_i=k_i \alpha^{-\lambda}.
\end{equation}
The only problem with (\ref{dfopid_ideal}) is that infinite number
of memory units are needed for its realization, and consequently
the computational cost is increased by increasing the time. In
other words, in practice the upper bound of sigmas in
(\ref{dfopid_ideal}) cannot be considered equal to infinity.
Restricting the number of memory units to $M$, the following
formula is proposed for the $M$th-order LDPID controller:
\begin{equation}\label{dfopid_prac}
C_d(z)=K_p+K_d \sum_{k=0}^M f_k(\mu)z^{-k}+ K_i
\frac{1+z^{-1}}{1-z^{-1}} \sum_{k=0}^M f_k(1-\lambda)z^{-k}.
\end{equation}
In the rest of this paper whenever we refer to the LDPID
controller, a system with transfer function (\ref{dfopid_prac}) is
under consideration. Figure \ref{fig_dfopid_prac} shows the block
diagram of the proposed LDPID controller where C/D and D/C stand
for the ideal continuous-to-discrete-time and
discrete-to-continuous-time converters, respectively
\cite{oppenheim2}. Obviously, in practice the C/D is realized
using a sample-and-hold (or an A/D converter) and the output of
the adder in Fig. \ref{fig_dfopid_prac} can directly be applied to
the process. For simulation in Matlab, the C/D is modeled with a
sample-and-hold, and the D/C is simply omitted. It will be shown
later that one advantage of the proposed LDPID controller is that
it can easily cope with the unwanted effects caused by the
sample-and-hold used in practice.

It is important to note that the proposed LDPID controller as
defined in (\ref{dfopid_prac}) actually has six parameters to
tune. Considering the fact that the continuous-time and
discrete-time FOPID controllers (as defined in
(\ref{fopid_transfer}) and (\ref{dfopid_prac}), respectively) have
different number of parameters to tune, and taking into account
the effect of sample-and-hold, it is evident that in practice the
parameters of these two controllers cannot simply be related
according to (\ref{kp}). In other words, the transfer function of
a certain LDPID controller is not, in general, obtained by
applying the Tustin transform to any FOPID controller. It
concludes that the parameters of the proposed LDPID controller
should be tuned directly. For this purpose, a certain value can be
assigned to $M$ and then the value of other parameters be
calculated such that a predetermined set of objectives is met. The
other possible approach is to consider $M$ as a tuning parameter
and then find the values of $K_p$, $K_d$, $K_i$, $\mu$, $\lambda$,
and $M$ such that the objectives under consideration are met. The
first approach is used in the rest of this paper.

\subsection{Two methods for tuning the proposed LDPID
controller}\label{sec_method} Two methods are developed in this
section for tuning the parameters of the LDPID controller defined
in (\ref{dfopid_prac}). The first method is the discrete-time
equivalent of the method proposed in \cite{monje2} for tuning the
parameters of FOPID controller. The second method is based on
minimization of a certain integral performance index by suitable
choice of the unknown parameters of controller. This approach is
similar to the method used in \cite{farshad1} for optimal tuning
the FOPID controllers. In the following we discuss on these two
methods briefly.

In both methods first we assign a suitable value to $M$. This
value should be chosen considering the limitations of the hardware
used to realize the controller. Evidently, increasing the value of
$M$ increases the computational cost and the memory usage, and
indeed puts a limitation on the minimum possible value for
sampling period of system. (Simulations performed by authors show
that the final results are not so sensitive to the special value
assigned to $M$ provided that other parameters are selected
properly.)

Note that according to (\ref{dfopid_prac}) the difference equation
relating $e[n]$ to $u[n]$ in Fig. \ref{fig_dfopid_prac} is as the
following:
\begin{multline}\label{dfopid_difference_eq}
u[n]=u[n-1]+K_p(e[n]-e[n-1])+ \sum_{k=0}^M K_d f_k(\mu)
\left(e[n-k]-e[n-k-1]\right)+\\  \sum_{k=0}^M K_i f_k(1-\lambda)
(e[n-k]-e[n-k-1]).
\end{multline}
In the above equation calculation of each sigma needs $M+1$
(floating point) multiplications and $M+1$ summations (assuming
that $K_d f_k(\mu)$ and $K_i f_k(1-\lambda)$ are calculated
beforehand and are known parameters). Hence, it can be easily
verified that calculation of $u[n]$ from
(\ref{dfopid_difference_eq}) totally needs  $2(M+1)+1=2M+3$
multiplications and $2(M+1)+4=2M+6$ summations at each sampling
period (note that (\ref{dfopid_difference_eq}) is not a minimal
description for (\ref{dfopid_prac}) and more effective
formulations can also be obtained \cite{oppenheim2}). Considering
the fact that any floating point multiplication is much more time
consuming than any summation, the computational cost of the
proposed controller can be approximated by the number of
multiplications, which is equal to $2M+6$ at each sampling period.
So, in order to determine the suitable value of $M$ first we
should determine the suitable sampling period of system and then
choose the value of $M$ such that the digital processor can
perform at least $2M+6$ floating point multiplications at each
sampling period.

After determining the value of $M$, the values of the remaining
five parameters of controller should be determined such that the
following five conditions are satisfied  simultaneously ($P(s)$ is
the process transfer function, which is located in a standard
unity feedback system in series with the LDPID controller):
\begin{itemize}
    \item The gain crossover frequency of the open-loop system,
    $\omega_c$, be equal to the desired value, that is the
    equality
    \begin{equation}\label{cond1}
\left|C_d(e^{j\omega_c}) P(j\omega_c)\right|=0\mathrm{dB},
    \end{equation}
    holds for the desired $\omega_c$, where $P(s)$ is the process
    transfer function.
    \item The phase margin of the feedback system, $\varphi_m$, be
    equal to the desired value, that is the equality
    \begin{equation}\label{cond2}
\arg\left\{C_d(e^{j\omega_c}) P(j\omega_c)\right\}=\pi +\varphi_m,
    \end{equation}
    holds for the desired $\varphi_m$.
    \item The feedback system exhibits a good robustness to variations in the
    gain of process, which can be achieved by satisfying the
    following equality
    \begin{equation}\label{cond3}
\frac{\mathrm{d}\left(C_d(e^{j\omega}) P(j\omega) \right)}
{\mathrm{d} \omega}\bigg|_{\omega = \omega_c}=0.
    \end{equation}
    \item The feedback system attenuates the high frequency noise,
    which is achieved by satisfying the inequality:
    \begin{equation}\label{cond4}
\left|\frac{C_d(e^{j\omega}) P(j\omega)} {1+C_d(e^{j\omega})
P(j\omega)}\right|\le A~ \mathrm{dB}\quad \omega\ge \omega_t
\mathrm{rad/s},
    \end{equation}
    where $A$ and $\omega_t$ are desired constants.
    \item The feedback system rejects the disturbance, which is achieved by satisfying the inequality:
    \begin{equation}\label{cond5}
\left|\frac{1} {1+C_d(e^{j\omega}) P(j\omega)}\right|\le B~
\mathrm{dB}\quad \omega\le \omega_s \mathrm{rad/s},
    \end{equation}
        where $B$ and $\omega_s$ are desired constants.
\end{itemize}
Similar to \cite{monje2}, in this paper (\ref{cond1}) is
considered as the main object of optimization and
(\ref{cond2})-(\ref{cond5}) are considered as the corresponding
constrains. More precisely, we have applied the genetic algorithm
to find the values of $K_p$, $K_d$, $K_i$, $\mu$, and $\lambda$ in
(\ref{dfopid_prac}) such that $\left|\left|C_d(e^{j\omega_c})
P(j\omega_c)\right|-1\right|$ is minimized and simultaneously the
equality constraints (\ref{cond2}) and (\ref{cond3}), and
inequality constraints (\ref{cond4}) and (\ref{cond5}) are
fulfilled. Clearly, it may happen that the above mentioned
optimization problem does not have any solution, but even
approximate solutions (which violate the constrains to some
extent) are useful in practice. Numerical simulations performed by
authors show that the genetic algorithm toolbox of Matlab can
effectively solve such a complicated constrained optimization
problem in a relatively short time (see Examples 1 of Section 3),
however the application of this tuning strategy is restricted to
relatively simple processes.

The second method that can be used for tuning the parameters of
LDPID controller is to assign a certain value to $M$ and then
calculate the values of $K_p$, $K_d$, $K_i$, $\mu$, and $\lambda$
in (\ref{dfopid_prac}) such that an integral performance index
(e.g., the IAE or ISE performance index corresponding to the
tracking error of step command) is minimized. In this method, the
genetic algorithm (or any other meta-heuristic optimization
algorithm such as PSO) can be used to search the five-dimensional
space to find the optimal solution. For this purpose, the
corresponding integral performance index should be considered as
the object of minimization and the equation obtained by equating
the number of unstable poles of the closed-loop system to zero can
be considered as the constraint of optimization. This equality
constraint is mandatory when the integral performance index is
evaluated in the frequency domain (using Parseval's theorem) since
in this case the stability issue is not involved in the cost
function, i.e., a set of parameters that minimize the integral
performance index in frequency domain may lead to an unstable
feedback system.

\subsection{Generalization}\label{sec_general}
The coefficients $f_k(\mu)$ and $f_k(1-\lambda)$ of the proposed
LDPID controller (\ref{dfopid_prac}) are calculated from
(\ref{weight}). The question that may arise at this point is: Can
one propose another reasonable method instead of (\ref{weight})
for calculation of these coefficients? To provide an answer for
this question first recall that any classical PID controller
(either of continuous or discrete type) combines a proportion,
derivative and integral of the error signal to generate the
control. It motivates us first to study the general behavior of
discrete-time derivative operators. Then based on these results we
can conclude that any alternating-sign function of $k$ provides us
with a natural definition for $f_k(\mu)$ to be used in
(\ref{dfopid_prac}). Finally, we briefly extend the results to
discrete-time integrators.

Suppose that we want to approximate the time derivative of the
function $e(t)$ only by using its samples $e_i=e(t_i)$
($i=1,\ldots,n$) where $T=t_i-t_{i-1}$. Mathematically, this task
can be performed by using one of the following formulas
\cite{hornbeck}:
\begin{equation}
\dot{e}(t_i)=\frac{e_{i+1}-e_i}{T}+O(T),
\end{equation}
\begin{equation}
\dot{e}(t_i)=\frac{e_i-e_{i-1}}{T}+O(T),
\end{equation}
\begin{equation}
\dot{e}(t_i)=\frac{e_{i+1}-e_{i-1}}{2T}+O(T^2),
\end{equation}
which are called forward, backward, and centered difference
approximations, respectively. In each case one can increase the
accuracy of the resulted time derivative by using larger number of
sample points. For example, backward difference approximations
with second and third order errors are obtained as the following
\cite{hornbeck}:
\begin{equation}
\dot{e}(t_i)=\frac{3e_i-4e_{i-1}+e_{i-2}}{2T}+O(T^2),
\end{equation}
\begin{equation}
\dot{e}(t_i)=\frac{11e_i-18e_{i-1}+9e_{i-2}-2e_{i-3}}{6T}+O(T^3).
\end{equation}
The general formula for approximating the time derivative of a
certain function with $n$th-order error using the backward
difference method is as the following \cite{khan}:
\begin{equation}\label{temp2}
\dot{e}(t_i)=\frac{1}{T}\sum_{k=-n}^0 g_{k,n}^{B,1}
e_{k+i}+O(T^n),
\end{equation}
where the coefficients $g_{k,n}^{B,1}$ are calculated from the
following iterative procedure:
\begin{equation}\label{temp3}
g_{0,n}^{B,1}=\sum_{j=1}^n (1/j),\quad g_{-1,n}^{B,1}=-n,
\end{equation}
\begin{equation}\label{temp1}
g_{-k,n}^{B,1}=-g_{-k+1,n}^{B,1}(k-1)(n-k+1)/k^2,\quad
k=2,\ldots,n.
\end{equation}
Very similar formulas can also be found in \cite{khan} for forward
and centered difference approximations. The intersection point of
all of these finite-difference formulas is that they approximate
the time derivative of a function by the weighted sum of its
samples, where these weights change sign decussately as the $k$ is
increased. For example, according to (\ref{temp1}) it is obvious
that the weights $g_{k,n}^{B,1}$ in (\ref{temp2}) are positive for
$k=0,-2,-4,\ldots$ and negative for $k=-1,-3,-5,\ldots$.

\begin{figure}[tb]
\begin{center}
\includegraphics[width=8.5cm]{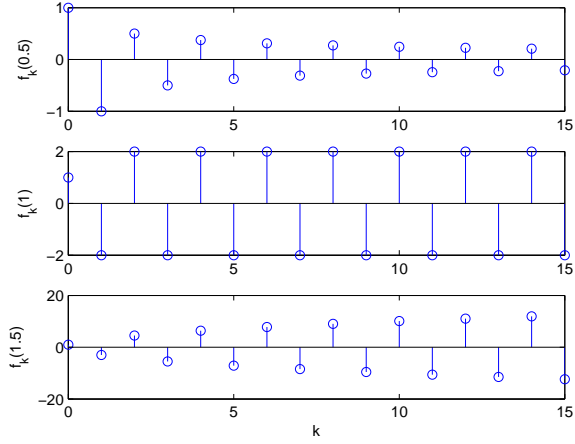}
\caption{The plot of $f_k(\mu)$ versus $k$ for three different
values of $\mu$.} \label{fig_fk1}
\end{center}
\end{figure}

Equation (\ref{temp2}) leads us to the fact that in digital
control system design one can use the following $n$th-order
difference equation to calculate the derivative of error signal
from its samples with an arbitrary precision:
\begin{equation}\label{temp4}
e_\mathrm{dot}[i]\approx \frac{1}{T}\sum_{k=0}^n g_{-k,n}^{B,1}
e[i-k],
\end{equation}
where $e[i]$ and $e_\mathrm{dot}[i]$ stand for the samples of
error signal and its derivative, respectively, and the
alternating-sign weights $g_{-k,n}^{B,1}$ are defined similar to
(\ref{temp3}) and (\ref{temp1}) (note that the same discussion
goes for forward and centered difference approximations, however,
the backward difference approximation is advantageous in the way
that it leads to a causal relation). Recall that according to
(\ref{dfopid_prac}) the proposed LDPID controller relates the
samples of error signal to its derivative through the following
difference equation
\begin{equation}\label{temp5}
e_\mathrm{dot}[i]=K_d \sum_{k=0}^M f_k(\mu)e[i-k],
\end{equation}
where the weights $f_k(\mu)$ are shown in Fig. \ref{fig_fk1} for
three different values of $\mu$. As it can be observed in this
figure the weights $f_k(\mu)$ in (\ref{temp5}) change sign
decussately similar to the weights $g_{-k,n}^{B,1}$ in
(\ref{temp4}). It means that the derivative term of the proposed
LDPID controller also has the property of subtracting many two
successive (positive-weighted) samples of the error signal and
then forming the cumulative sum of results. However, the
difference equation (\ref{temp5}) is, compared to (\ref{temp4}),
advantageous in the way that one can adjust the amplitude of
weights simply by changing the value of $\mu$ (note that in case
of using (\ref{temp4}) we have no control on the amplitude of the
weights $g_{-k,n}^{B,1}$ for the given $n$). For this reason, it
is not surprising if the derivative action of the LDPID controller
as given in (\ref{temp5}) reduces to a classical discrete-time
differentiator for some $\mu$ and $M$. For example, assuming
$\mu=0.5$ and $M=1$ equation (\ref{temp5}) yields
\begin{equation}
e_\mathrm{dot}[i]=K_d \sum_{k=0}^1 f_k(0.5)e[i-k]=K_d
(e[i]-e[i-1]),
\end{equation}
which is the classical backward difference approximation for
derivative operator.

The above discussion motivates us to define, in general, the
derivative action of discrete-time LDPID controllers through the
difference equation (\ref{temp5}) where $f_k(\mu)$ ($\mu>0$) can
be considered equal to any alternating-sign function of $k$ (which
is not necessarily calculated from (\ref{weight})). This
definition for discrete-time derivative operator is consistent
with the one used in Section \ref{sec_ldpid} and the classical
higher-order one presented in (\ref{temp4}). However, application
of, e.g., the backward difference approximation in (\ref{tus1})
leads to a discrete-time approximation for derivative operator
that is not consistent with the above generalized definition and a
more general definition is needed to cover this case (for example,
by letting $f_k$ to take the positive and negative values
arbitrarily).

The integral term of the proposed LDPID controller can be
generalized very similar to the above mentioned approach. More
precisely, the generalized discrete-time integrator can be defined
as an operator that generates the positive-weighted sum of the
current and past samples of error signal (it is not difficult to
show that the last term in the right hand side of
(\ref{dfopid_prac}) generates the positive-weighted sum of the
samples of error signal for $\lambda>0$).

\section{Simulated and Practical Examples}
Four numerical examples and a real-world application are studied
in this section to verify the effectiveness of the proposed LDPID
controller and in support of its superiority over the classical
PID controllers.\\

\begin{figure}[tb]
\begin{center}
\includegraphics[width=8.5cm]{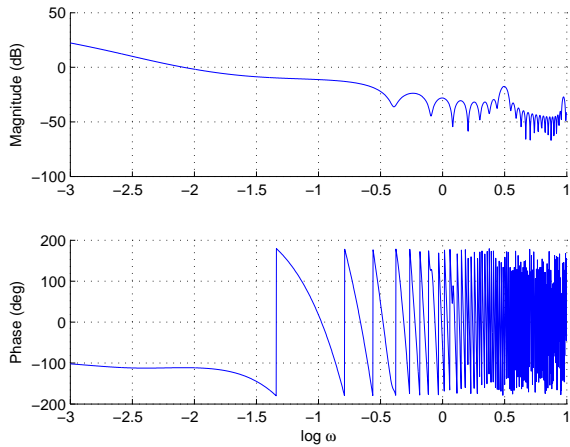}
\caption{Bode plot of $C_d(e^{sT})G(s)$, corresponding to Example
1.} \label{fig_ex2_1}
\end{center}
\end{figure}

\begin{figure}[tb]
\begin{center}
\includegraphics[width=8.5cm]{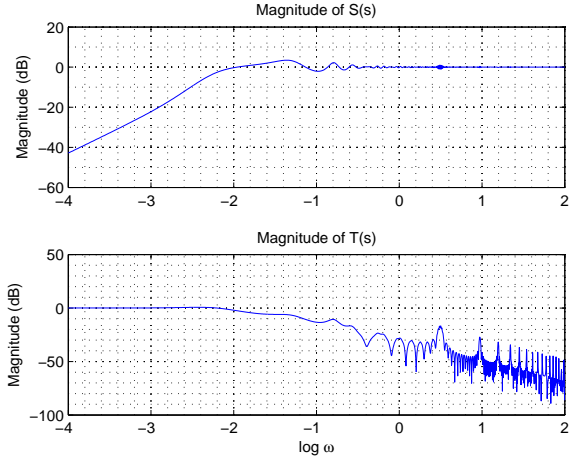}
\caption{Bode magnitude plot of $S(s)$ and $T(s)$, corresponding
to Example 1.} \label{fig_ex2_3}
\end{center}
\end{figure}

\begin{figure}[tb]
\begin{center}
\includegraphics[width=8.5cm]{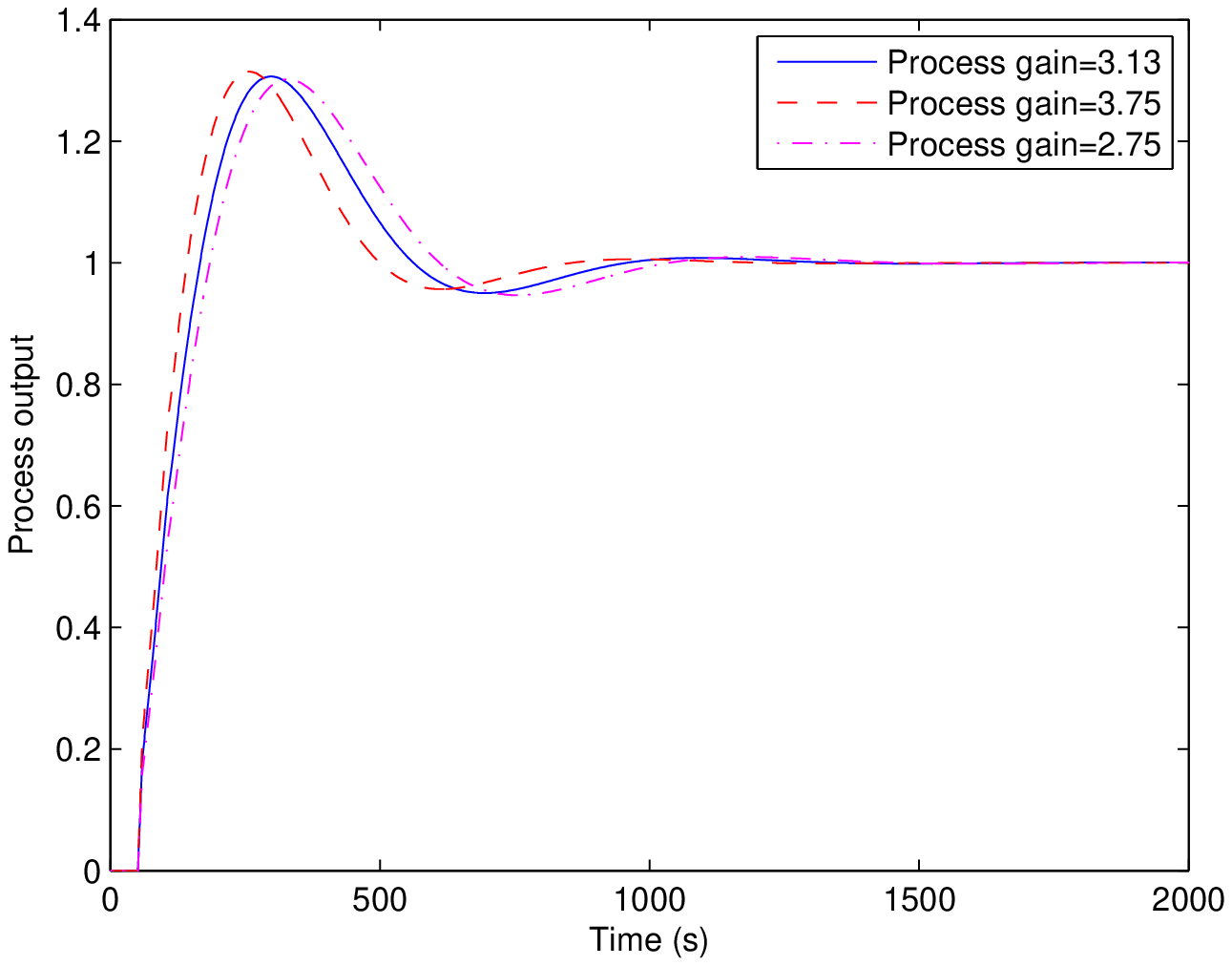}
\caption{Unit step response of the feedback system when the LDPID
controller is applied, corresponding to Example 1.}
\label{fig_ex2_2}
\end{center}
\end{figure}

\textbf{Example 1: Application to a varying-gain process.} The
problem under consideration in this example is adopted form
\cite{monje2}. Consider a first-order plus time-delay (FOPTD)
process with the following transfer function:
\begin{equation}
P(s)=\frac{Ke^{-50s}}{433.33s+1},
\end{equation}
where the nominal value of $K$ is equal to $3.13$ and it may vary
between $2.75$ and $3.75$. Our aim here is to design a LDPID
controller such that the following objectives are achieved
simultaneously:
\begin{itemize}
    \item $\omega_{c}=0.008$ rad/s,
    \item $\varphi_m=60^\circ$,
    \item $|S(j\omega)|\le -20$ dB, $\forall\omega\le \omega_s
    =0.001$
    rad/s,
    \item $|T(j\omega)|\le -20$ dB, $\forall\omega\ge \omega_t =10$
    rad/s,
    \item robustness to variations in the gain of the process.
\end{itemize}
Following the design procedure presented in Section 2 the
following LDPID controller is obtained (assuming $M=15$) using the
genetic algorithm toolbox of Matlab:
\begin{multline}\label{c1_ex2}
C_d(z)=3.059+0.384\sum_{k=0}^{15} f_k(1.228) z^{-k} +0.059
\frac{1+z^{-1}}{1-z^{-1}} \sum_{k=0}^{15} f_k(0.45) z^{-k}.
\end{multline}
Note that the above controller is actually a very good suboptimal
solution which approximately satisfies all requirements of the
problem. Figure \ref{fig_ex2_1} shows the Bode plot of
$C_d(e^{sT})G(s)$. In this figure the gain crossover frequency is
approximately equal to $0.008$ rad/s and the phase plot is almost
flat around this frequency, as it is expected (note that
$\log_{10} 0.008\approx -2.1$). Moreover, the phase margin is
slightly larger than $60^\circ$. Figure \ref{fig_ex2_3} shows the
Bode magnitude plots of the sensitivity and complementary
sensitivity function. Using this figure it can be easily verified
that $S(s)$ and $T(s)$ satisfy the requirements of problem. Figure
\ref{fig_ex2_2} shows the unit step response of the corresponding
feedback system for three different values of $K$. As it can be
observed, the overshoot of response is almost constant for
different values of $K$, which is the direct consequent of the
flatness of the Bode phase plot of open-loop transfer function at
frequencies around $\omega_{c}$.\\

\begin{figure}[tb]
\begin{center}
\includegraphics[width=8.5cm]{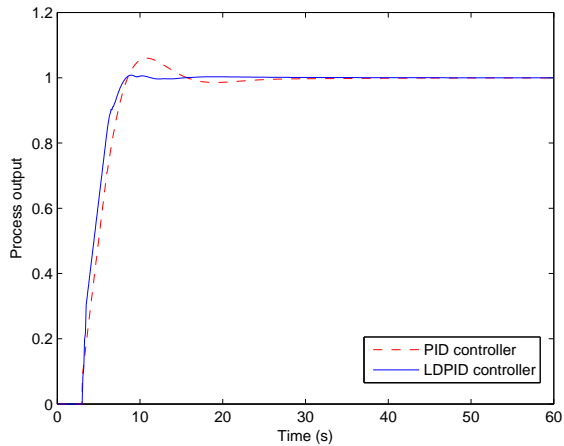}
\caption{Unit step response of the closed-loop system,
corresponding to Example 2.} \label{fig_ex1_2}
\end{center}
\end{figure}

\textbf{Example 2: Application of the optimal LDPID controller to
a FOPTD process.} The aim of this example is to show that
designing a classical PID controller for a FOPTD process and then
discretizing it (using, e.g., the Tustin transform with
prewarping) may lead to an unstable feedback system, while direct
tuning the proposed LDPID controller can remove this difficulty.
For this purpose consider a FOPTD process with the following
transfer function
\begin{equation}\label{p_ex1}
P(s)=\frac{2e^{-3s}}{1+10s}.
\end{equation}
This process can effectively be controlled with the following PID
controller, which is obtained by trial and error:
\begin{equation}\label{c_ex1}
C(s)=1.1+\frac{0.1}{s}+0.4s.
\end{equation}
Figure \ref{fig_ex1_2} shows the unit step response of the
feedback system in this case. As it can be observed, the command
following is quite satisfactory.

\begin{figure}[tb]
\begin{center}
\includegraphics[width=8.5cm]{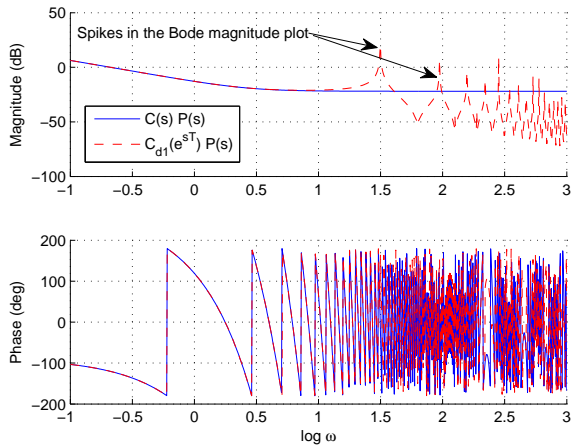}
\caption{Bode plots of $C(s)P(s)$ and $C_{d1}(e^{sT})P(s)$,
corresponding to Example 2.} \label{fig_ex1_1}
\end{center}
\end{figure}

Now let us examine the performance of this feedback system when
controller (\ref{c_ex1}) is discretized and then applied. In this
example the gain crossover frequency of the open-loop system is
$\omega_c\approx 0.21$ rad/s and the phase margin is about
$60^\circ$. Applying prewarped Tustin method to (\ref{c_ex1})
(assuming that the sampling period is equal to $T=0.1$ s) leads to
the following discrete-time controller:
\begin{equation}\label{c2_ex1}
C_{d1}(z)=1.1+0.005 \frac{1+z^{-1}}{1-z^{-1}} + 8
\frac{1-z^{-1}}{1+z^{-1}}.
\end{equation}
Figure \ref{fig_ex1_1} shows the Bode plots of $C(s)P(s)$ and
$C_{d1}(e^{sT})P(s)$ (the latter corresponds to the open-loop
transfer function when the discrete-time controller is applied in
series with ideal C/D and D/C converters). This figure clearly
shows that the closed-loop system with discrete-time controller is
unstable, while we observed that the original continuous-time
feedback system is stable with a satisfactory phase margin (recall
that any spike in the Bode magnitude plot of $C_{d1}(e^{sT})P(s)$
in Fig. \ref{fig_ex1_1} corresponds to an encirclement around -1
in the corresponding Nyquist plot). It should be emphasized that
the dashed curve in Fig. \ref{fig_ex1_1} is obtained assuming that
ideal C/D and D/C converters are available. Of course, applying a
practical A/D will make the open-loop system (with discrete-time
controller) even more phase lag. One can also perform a
time-domain simulation to justify the fact that the discretized
controller (\ref{c2_ex1}) leads to an unstable feedback system.
For this purpose the error signal can be discretized using the
Zero-Order Hold box in Simulink, and then the resulted signal is
considered as the input of (\ref{c2_ex1}), and the corresponding
output is entered to the input of plant. It is worth to mention
that in this example one can hardly arrive at a stable feedback
system only by changing to sampling period used for discretization
of the controller.

Assuming $M=5$ in (\ref{dfopid_prac}), the (suboptimal) LDPID
controller which minimizes the IAE performance index
(corresponding to the tracking error of the unit step command) is
obtained as the following:
\begin{multline}\label{c3_ex1}
C_{d2}(z)=2.8+1.5\sum_{k=0}^{5} f_k(1.03) z^{-k} +0.004
\frac{1+z^{-1}}{1-z^{-1}} \sum_{k=0}^{5} f_k(-0.1) z^{-k}.
\end{multline}
The main advantage of this controller over the one given in
(\ref{c_ex1}) is that the unwanted effects caused by using
zero-order hold are also taken into account during the controller
design and consequently it is ready to be realized using
microprocessors. Figure \ref{fig_ex1_2} also shows the unit step
response of the closed-loop system when the LDPID controller
(\ref{c3_ex1}) is applied. As it can be observed, the LDPID
controller results in
a satisfactory transient response.\\

\textbf{Example 3: Application to a non-minimum phase process with
oscillatory poles.} The following transfer function appears in the
one-link flexible arm robot \cite{chen2}:
\begin{equation}\label{p_ex3}
P(s)=\frac{-4.906s^2-0.5884s+335.17}{s^4+0.55437s^3+139.6s^2+27.91s}.
\end{equation}
$P(s)$ has a non-minimum phase zero at $8.2057$ and four poles at
0, $0.2$, $0.1772 \pm j11.8109$. Controlling a system with
transfer function (\ref{p_ex3}) is a relatively difficult task
since trivial PID controllers often do not lead to satisfactory
results when the process has both the non-minimum phase zero and
complex conjugate poles with a very small damping ratio (here we
have $\zeta= 0.0150$) \cite{astrom2}. In the following we try the
proposed LDPID controller.

\begin{figure}[tb]
\begin{center}
\includegraphics[width=8.5cm]{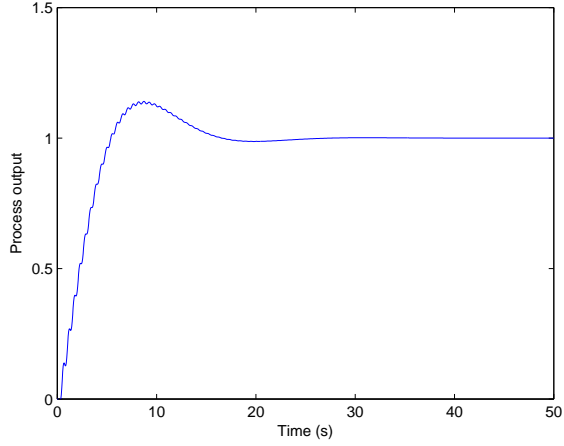}
\caption{Unit step response of the closed-loop system,
corresponding to Example 3.} \label{fig_ex3_1}
\end{center}
\end{figure}

Since in this example the process itself has a pole at the origin,
a LDPD controller is sufficient for the tracking of step command
without steady-state error (the LDPD controller is obtained by
omitting the integral term of (\ref{dfopid_prac})). In order to
design the LDPD controller first we arbitrarily assume $M=5$ and
then we obtain the parameters of the LDPD such that the IAE
performance index (corresponding to the tracking error of step
command) is minimized. Next, we slightly modify the parameters of
the resulted controller by trial and error such that the
closed-loop system exhibits the desired step response. This
approach leads to the following LDPD controller:
\begin{equation}
C_d(z)=0.3+0.5\sum_{k=0}^5 f_k(0.8)z^{-k}.
\end{equation}
Unit step response of the corresponding closed-loop system is
shown in Fig. \ref{fig_ex3_1}. In this figure the rise time, the
settling time and overshoot of the response are approximately
equal to 5.4 s, 15 s, and $14\%$, respectively. Numerous numerical
simulations performed by authors show that a response with this
characteristics cannot be achieved by applying any PD or PID type
controller (i.e., either the rise time or the settling time or the
overshoot of the response obtained by using a PID controller will
be larger than the corresponding one obtained by using the
proposed LDPD controller). Interesting point is that the proposed
LDPD controller has, similar to PID controllers, only three
parameters to tune with the difference that it applies few more
memory units.\\

\begin{figure}[tb]
\begin{center}
\includegraphics[width=8.5cm]{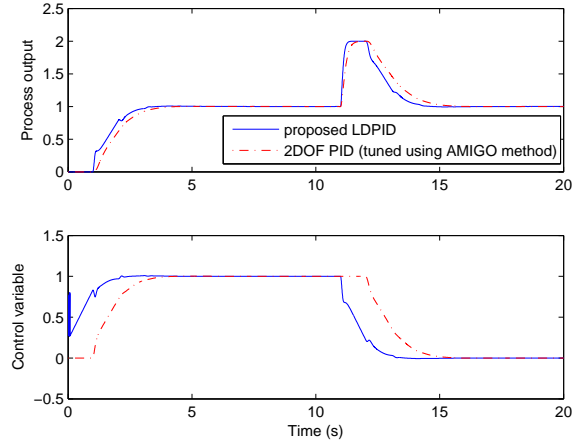}
\caption{Process output and control variable, corresponding to
Example 4.} \label{fig_ex4}
\end{center}
\end{figure}

\textbf{Example 4: Comparing LDPID with a 2DOF PID tuned by AMIGO
method.} \AA str\"{o}m and H\"{a}gglund \cite{astrom1} proposed
the so-called AMIGO method for tuning the two degree-of-freedom
(2DOF) PID controllers, which have six parameters to tune.
Application of this method to a process with transfer function
$P(s)=e^{-s}/(1+0.05s)^2$ leads to a 2DOF PID controller with the
following input-output relation:
\begin{multline}\label{amigo}
u(t)=0.242(y_{sp}-y_f)+0.515\int_{0}^t(y_{sp}(\tau)-y_f(\tau))
\mathrm{d}\tau-0.032\frac{\mathrm{d}y_f(t)}{\mathrm{d}t},
\end{multline}
where $y_{sp}$ is the set point, $Y_f(s)=G_f(s)Y(s)$ the filtered
process variable, $u(t)$ the control, $G_f(s)=1/(1+0.1s)^2$, and
$y(t)$ is the process output. On the other hand, minimization of
the IAE performance index (more precisely, the integral of the
absolute error when the unit step command and the unit step
disturbance are applied at $t=0$ and $t=10$, respectively in the
standard feedback connection) leads to the following (suboptimal)
LDPID controller:
\begin{multline}\label{cd_ex4}
C_d(z)=0.5+0.15\sum_{k=0}^{15} f_k(1.15) z^{-k} +7 \times 10^{-4}
\frac{1+z^{-1}}{1-z^{-1}} \sum_{k=0}^{15} f_k(-0.2) z^{-k}.
\end{multline}
Figure \ref{fig_ex4} shows the process output and the
corresponding control variable when the PID controllers
(\ref{amigo}) and (\ref{cd_ex4}) are applied. As it can be
observed in this figure both controllers lead to almost the same
step response, and both of them apply almost the same control
effort. It concludes that the proposed controller can compete the
2DOF PID controller tuned by using the AMIGO method. However, the
proposed controller is advantageous in the way that it is in the
discrete-time form and ready for realization.\\

\textbf{Example 5: Experimental results.} In this example we study
the application of the proposed LDPID controller for temperature
control of an industrial heating box. In this system a DC voltage
is applied to the heating element wire of the box and a PWM with
adjustable duty cycle is used to control the on-off time of this
wire. The DC voltage applied to wire does not have a certain
value. In fact, an industrial rectifier provides a very high DC
voltage which is applied to the series connection of the wires of
many boxes of this type. Hence, the DC voltage across the wire of
each box (when it is on) strictly depends on the number of wires
connected in series (typically, in practice few hundred boxes are
used in series).

For many reasons this process constitutes a difficult control
problem. Firstly, the DC voltage applied to each wire (which, of
course, affects the DC gain of process) has an uncertain value. In
fact, in practice it is observed that this voltage may even be
subjected to 200\% changes from the nominal value. It is also
observed that changing this voltage also changes the time constant
of process. Secondly, the time constant of the process also
depends on the type and weight of the material embedded in the box
(to be warmed). Thirdly, taking into account the function of PWM
it is obvious that the system is actually nonlinear. More
precisely, the negative control signals generated by controller
are simply neglected by PWM. Note also that applying a kind of
anti-windup technique is also mandatory for digital realization of
integrators.

At different working conditions each box can (very approximately)
be modelled by the following first-order uncertain transfer
function
\begin{equation}
P(s)=\frac{K}{1+sT},\quad 14<K<34, \quad 380<T<570.
\end{equation}
(In the above transfer function the input is duty cycle and the
output is temperature in $^\circ$C). The nominal process model is
also considered as $P(s)=32/(1+425s)$ (note that the probability
of occurring different uncertainties is not the same). Our aim
here is to design a controller which leads to $\omega_c=1$ rad/s,
$\varphi_m\ge 75^\circ$, $\omega_t=10$ rad/s, $\omega_s=0.1$
rad/s, and $A=B=- 20$ dB (see Section \ref{sec_method}). Moreover,
according to the high uncertainty in the process model it is
highly desired that the open-loop phase plot be as flat as
possible at frequencies around $\omega_c$.

\begin{figure}[tb]
\begin{center}
\includegraphics[width=8.5cm]{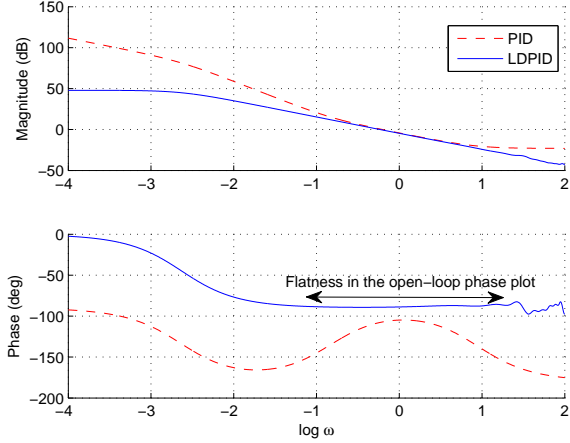}
\caption{Bode plots of $C(s)P(s)$ and $C_d(e^{sT})P(s)$,
corresponding to Example 5.}\label{fig_ex50}
\end{center}
\end{figure}

Following the procedure presented in Section \ref{sec_method} the
genetic algorithm leads to the following PID and LDPID controllers
($T=0.1$s):
\begin{equation}
C(s)=7.937+\frac{1.187}{s}-0.935s,
\end{equation}
\begin{multline}
C_d(z)=7.109+0.711\sum_{k=0}^5
f_k(0.077)z^{-k}+0.750\frac{1+z^{-1}}{1-z^{-1}}\sum_{k=0}^5
f_k(0.415)z^{-k}.
\end{multline}
Figure \ref{fig_ex50} shows the Bode phase and magnitude plots of
$C(s)P(s)$ and $C_d(e^{sT})P(s)$ when the nominal process model is
considered. In this example PID and LDPID controllers lead to
phase margins equal to $73.3^\circ$ and $91^\circ$, respectively.
As it can be observed in this figure the LDPID controller has
perfectly satisfied the design requests. Especially, unlike the
PID controller, the proposed LDPID has led to a very flat curve in
the phase plot (for more than 2 decades), which is highly desired
in dealing with the uncertain process under consideration. Using
trivial simulations it can be easily verified that both
controllers lead to very similar and satisfactory time-domain
responses when the nominal linear process model is considered.

\begin{figure}[tb]
\begin{center}
\includegraphics[width=8cm]{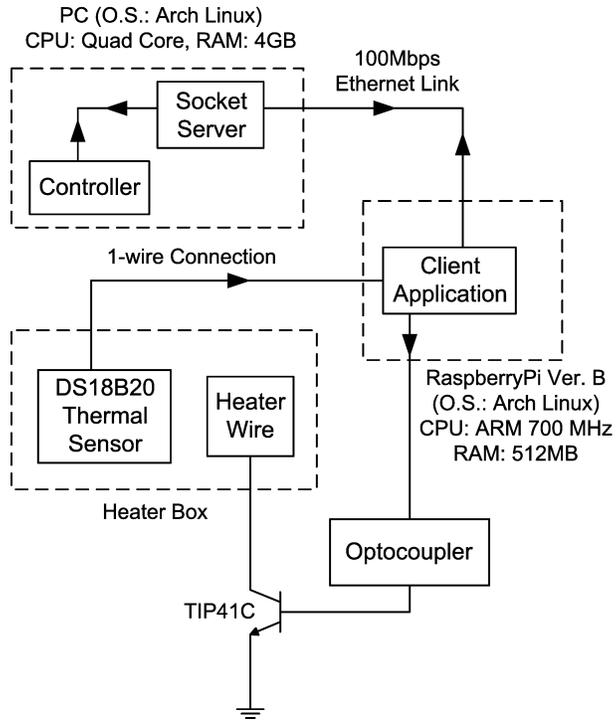}
\caption{The digital system used to realize the controllers of
Example 5.}\label{fig_ex53}
\end{center}
\end{figure}

Both of the above mentioned controllers are realized using the
digital system shown in Figure \ref{fig_ex53} (the classical PID
is discretized using the Tustin method with pre-wrapping
beforehand). Here it is worth to mention that since the 1-wire
output of DS18B20 digital sensor cannot directly be connected to
PC, a kind of transducer is needed. The RaspberryPi in Fig.
\ref{fig_ex53} is used for this purpose. The input of optocoupler
is connected to the software-generated PWM with frequency 50Hz. It
was observed that in dealing with LDPID controller calculation of
each control signal takes about 2ms in practice.

\begin{figure}[tb]
\begin{center}
\includegraphics[width=8.5cm]{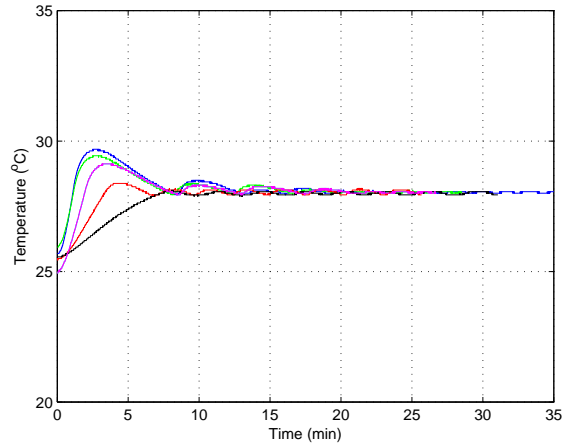}
\caption{Step responses of the practical closed-loop system when
LDPID is applied, corresponding to Example 5.}\label{fig_ex51}
\end{center}
\end{figure}

\begin{figure}[tb]
\begin{center}
\includegraphics[width=8.5cm]{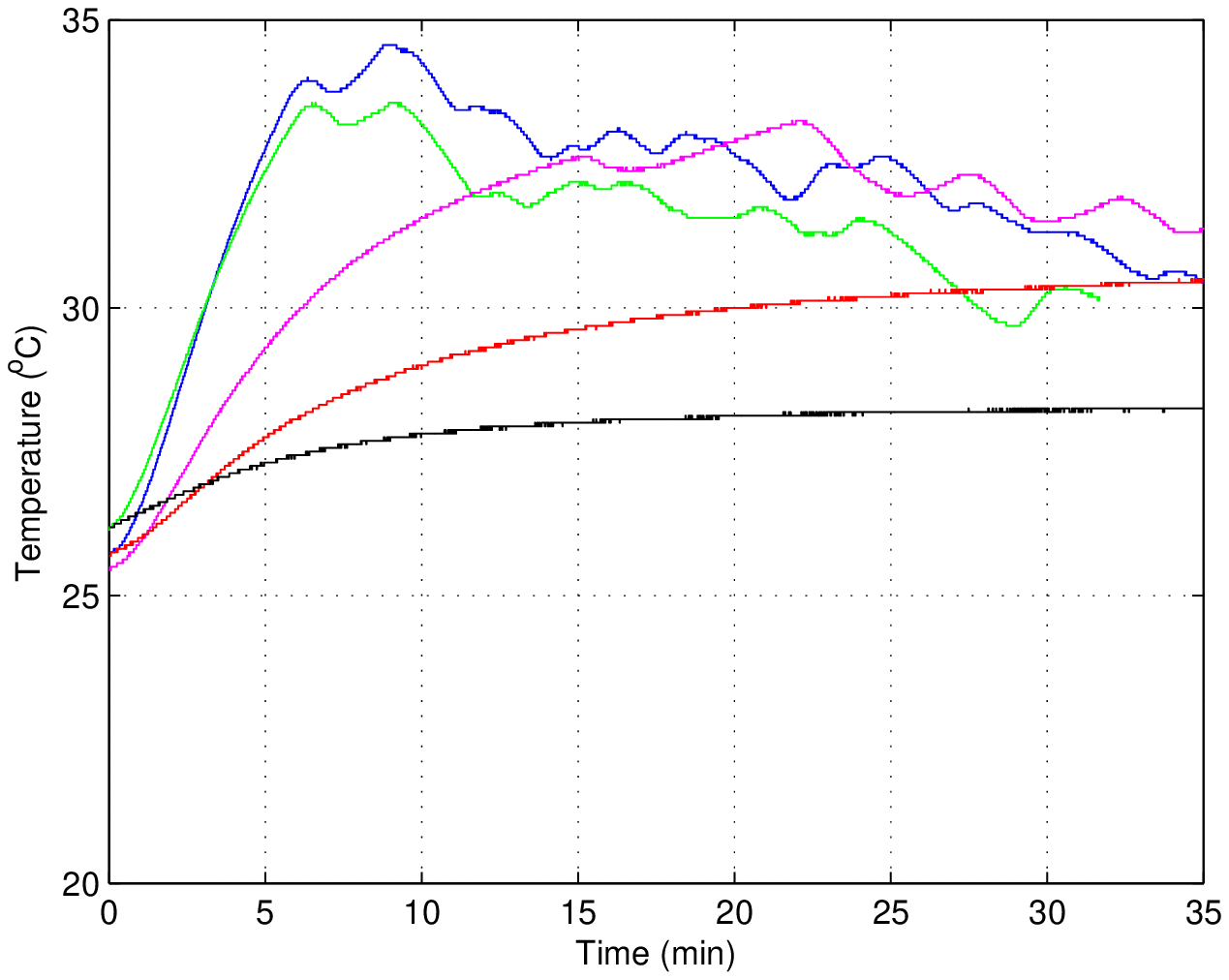}
\caption{Step responses of the practical closed-loop system when
classical PID is applied, corresponding to Example
5.}\label{fig_ex52}
\end{center}
\end{figure}

The time-domain responses of the practical closed-loop system when
PID and LDPID controllers are applied (assuming that the reference
temperature is equal to $28^\circ$C) are shown in Figs.
\ref{fig_ex51} and \ref{fig_ex52}, respectively in five different
conditions. Note that the curves with similar colors in Figs.
\ref{fig_ex51} and \ref{fig_ex52} are obtained under exactly the
same conditions in practice (i.e., exactly the same DC voltages
across the heating wire, the same materials in the box, etc.).
Note also that the system itself has an initial condition and the
vertical axis in Figs. \ref{fig_ex51} and \ref{fig_ex52} begins
from $20^\circ$C. Figure \ref{fig_ex51} clearly shows the
superiority of the proposed LDPID controller. In fact, the maximum
overshoot (in the worst case) caused by PID and LDPID controllers
is equal to 23.43\% and 6.1\%, respectively. Moreover,
fluctuations in the response are settled down much faster when
LDPID is applied.

It should be emphasized that since the PWM cannot generate
negative voltages (i.e., the process only has a heater but not a
cooler) it is observed in Fig. \ref{fig_ex51} that the temperature
is increased and decreased with two different time-constants. Note
also that according to the thermal capacity of heater, the
temperature in the box keeps increasing even after turning off the
heater. Consequently, since the derivative term of classical PID
applies much larger controls compared to LDPID, it leads to larger
overshoots and settling times in the response as it is observed in
Fig. \ref{fig_ex52}.


\section{Conclusions}
In this paper we proposed a new formulation for PID controllers,
which can be considered as the discrete-time equivalent of
fractional-order PIDs. Experimental and numerical examples were
also presented which verified the fact that the proposed
controller is capable of solving some difficult control problems
and has some advantages to classical PID controllers. We also
developed two methods for tuning the parameters of the proposed
controller.

\end{document}